\RequirePackage{ifpdf}
\ifpdf 
\documentclass[pdftex]{sigma}
\else
\documentclass{sigma}
\fi

\begin{document}
\allowdisplaybreaks

\renewcommand{\PaperNumber}{055}

\FirstPageHeading

\ShortArticleName{On the Existence of Conf\/igurations of
Subspaces in a Hilbert Space with Fixed Angles}

\ArticleName{On the Existence of Conf\/igurations of Subspaces\\
in a Hilbert Space with Fixed Angles}

\Author{Natasha D. POPOVA and Yurii S. SAMO\v{I}LENKO}

\AuthorNameForHeading{N.D. Popova and Yu.S. Samo\v\i{}lenko}

\Address{Institute of Mathematics, National Academy of  Sciences
of Ukraine, 3 Tereshchenkivs'ka Str., Kyiv-4, 01601 Ukraine}
\Email{\href{mailto:popova_n@yahoo.com}{popova{\_}n@yahoo.com},
\href{mailto:yurii_sam@imath.kiev.ua}{yurii\_sam@imath.kiev.ua}}

\ArticleDates{Received December 01, 2005, in f\/inal form April
30, 2006; Published online May 29, 2006}

\Abstract{For a class of $*$-algebras, where $*$-algebra
$A_{\Gamma,\tau}$ is generated by projections associated with
vertices of graph $\Gamma$ and depends on a parameter $\tau$ $(0 <
\tau \leq 1)$, we study the sets $\Sigma_\Gamma$ of values of
$\tau$ such that the algebras $A_{\Gamma,\tau}$ have nontrivial
$*$-representations, by using the theory of spectra of graphs. In
other words, we study such values of $\tau$ that the corresponding
conf\/igurations of subspaces in a Hilbert space exist.}

\Keywords{representations of $*$-algebras; Temperley--Lieb
algebras}

\Classification{16G99; 20C08}

\section{Introduction}
A number of papers is devoted to the study of $n$-tuples of
subspaces in a Hilbert space. The interest to this problem arose
in many respects due to its applications to problems of
mathematical physics (see, e.g., \cite{Evans&Kawahigashi} and the
bibliography therein).

Unitary description of $n$-tuples of subspaces is $*$-wild problem
when $n\geq 3$ (see \cite{Ostrovkyi&Samoilenko}).

In the present paper we study the conf\/igurations of subspaces
$H_i$ $(i=1,\ldots ,n)$ associated with the vertices of graph
$\Gamma$, where an angle between any two of subspaces is f\/ixed
(see Section~\ref{Subsp_in_Hilbert_spase}). It is convenient to
consider such conf\/igurations of subspaces as $*$-representations
of algebras generated by projections with relations of
Temperley--Lieb type (see \cite{Fan&Green,Temperley-Lieb,Wenzl}).
For tree~$\Gamma$ the set~$\Sigma_\Gamma$ (of those values of an
``angle'' $\tau$ where the corresponding conf\/igurations exist)
is described in Section~\ref{OnSetSigma}. This result is obtained
by using the theory of graph spectra (needed notions and results
are given in Section~\ref{about_graphs}). For graphs containing
cycles the situation is more complicated (see Remark~\ref{rem2}).

\section{Necessary facts from theory of spectra of graphs} \label{about_graphs}

We give some facts necessary for the exposition below, which can
be found in \cite{cvetkovic}. Let $\Gamma$ be a~f\/inite
undirected graph without multiple edges and loops. The
\textit{adjacency matrix} of a graph~$\Gamma$, with vertex set $\{
1,\ldots , n\}$, is $n\times n$ matrix
$A_{\Gamma}=\|a_{i,j}\|_{i,j=1}^n$ with $a_{i,j}=1$ if there is an
edge between~$i$ and $j$, and $a_{i,j}=0$; otherwise $a_{i,i}=0$
$\forall \, i$. The eigenvalues of $A_{\Gamma}$ and the spectrum
of~$A_{\Gamma}$ are also called the \textit{eigenvalues} and the
\textit{spectrum} of a graph~$\Gamma$, respectively. The
eigenvalues of $\Gamma$ are denoted by $\lambda_1,\ldots
,\lambda_n$; they are real because $A_{\Gamma}$ is symmetric. We
assume that $r=\lambda_1\geq\lambda_2\geq\cdots\geq\lambda_n=q.$
The largest eigenvalue $r=\lambda_1$ is called the \textit{index}
of a graph~$\Gamma.$

\begin{proposition}\label{inequalities}
$1.$  If a graph $\Gamma$ contains at least one edge then $1\leq
r\leq n-1$, $-r\leq q\leq -1$, and $r=-q$ if and only if a
component of $\Gamma$ with greatest index is a bipartite graph.

$2.$ If $\Gamma$ is a connected graph then $2\cos
\frac{\pi}{n+1}\leq r.$
\end{proposition}

\begin{remark}\label{q=-r}
If $\Gamma$ is a tree then $r=-q$, as the tree is a bipartite
graph.
\end{remark}

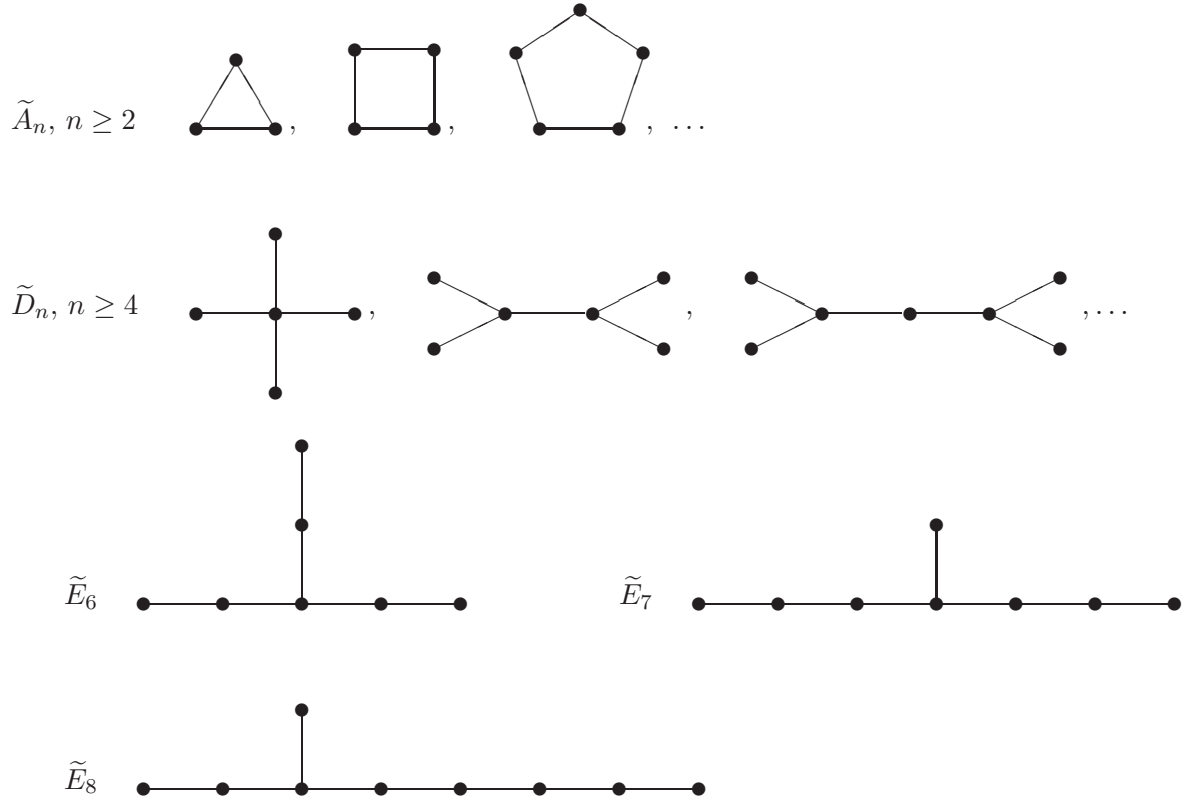
\begin{figure}[t]
\hspace*{3mm}\begin{picture}(440,305)
\put(20,0){$\widetilde{E}_8$} \put(50,0){\circle*{5}}
\put(50,0){\line(1,0){30}} \put(80,0){\circle*{5}}
\put(80,0){\line(1,0){30}} \put(110,0){\circle*{5}}
\put(110,0){\line(0,1){30}} \put(110,30){\circle*{5}}
\put(110,0){\line(1,0){30}} \put(140,0){\circle*{5}}
\put(140,0){\line(1,0){30}} \put(170,0){\circle*{5}}
\put(170,0){\line(1,0){30}} \put(200,0){\circle*{5}}
\put(200,0){\line(1,0){30}} \put(230,0){\circle*{5}}
\put(230,0){\line(1,0){30}} \put(260,0){\circle*{5}}
\put(20,70){$\widetilde{E}_6$} \put(50,70){\circle*{5}}
\put(50,70){\line(1,0){30}} \put(80,70){\circle*{5}}
\put(80,70){\line(1,0){30}} \put(110,70){\circle*{5}}
\put(110,70){\line(0,1){30}} \put(110,100){\circle*{5}}
\put(110,100){\line(0,1){30}} \put(110,130){\circle*{5}}
\put(110,70){\line(1,0){30}} \put(140,70){\circle*{5}}
\put(140,70){\line(1,0){30}} \put(170,70){\circle*{5}}
\put(230,70){$\widetilde{E}_7$} \put(260,70){\circle*{5}}
\put(260,70){\line(1,0){30}} \put(290,70){\circle*{5}}
\put(290,70){\line(1,0){30}} \put(320,70){\circle*{5}}
\put(320,70){\line(1,0){30}} \put(350,70){\circle*{5}}
\put(350,70){\line(0,1){30}} \put(350,100){\circle*{5}}
\put(350,70){\line(1,0){30}} \put(380,70){\circle*{5}}
\put(380,70){\line(1,0){30}} \put(410,70){\circle*{5}}
\put(410,70){\line(1,0){30}} \put(440,70){\circle*{5}}
\put(0,180){$\widetilde{D}_n,\, n\geq 4$}
\put(70,180){\circle*{5}} \put(70,180){\line(1,0){30}}
\put(100,180){\circle*{5}} \put(100,180){\line(0,1){30}}
\put(100,210){\circle*{5}} \put(100,180){\line(0,-1){30}}
\put(100,150){\circle*{5}} \put(100,180){\line(1,0){30}}
\put(130,180){\circle*{5}} \put(135,180){,} 
\put(160,193.41){\circle*{5}} \put(160,193.41){\line(2,-1){26.8}}
\put(160,166.59){\circle*{5}} \put(160,166.59){\line(2,1){26.8}}
\put(186.8,180){\circle*{5}} \put(186.8,180){\line(1,0){30}}
\put(220,180){\circle*{5}} \put(220,180){\line(2,-1){26.8}}
\put(246.8,193.41){\circle*{5}} \put(220,180){\line(2,1){26.8}}
\put(246.8,166.59){\circle*{5}}
\put(255,180){,}  
\put(280,193.41){\circle*{5}} \put(280,193.41){\line(2,-1){26.8}}
\put(280,166.59){\circle*{5}} \put(280,166.59){\line(2,1){26.8}}
\put(306.8,180){\circle*{5}} \put(306.8,180){\line(1,0){30}}
\put(340,180){\circle*{5}} \put(340,180){\line(1,0){30}}
\put(370,180){\circle*{5}} \put(370,180){\line(2,-1){26.8}}
\put(396.8,193.41){\circle*{5}} \put(370,180){\line(2,1){26.8}}
\put(396.8,166.59){\circle*{5}} \put(405,180){,}
\put(410,180){$\ldots$}
\put(0,250){$\widetilde{A}_n, \, n\geq 2$}
\put(70,250){\circle*{5}} \put(70,250){\line(1,0){30}}
\put(100,250){\circle*{5}} \put(70,250){\line(3,5){15}}
\put(100,250){\line(-3,5){15}} \put(85,276){\circle*{5}}
\put(105,250){,}
\put(130,250){\circle*{5}} \put(130,250){\line(1,0){30}}
\put(160,250){\circle*{5}} \put(130,250){\line(0,1){30}}
\put(130,280){\circle*{5}} \put(130,280){\line(1,0){30}}
\put(160,280){\circle*{5}} \put(160,280){\line(0,-1){30}}
\put(165,250){,}
\put(200,250){\circle*{5}} \put(200,250){\line(1,0){30}}
\put(230,250){\circle*{5}} \put(200,250){\line(-1,3){9}}
\put(191,278.5){\circle*{5}} \put(191,278.5){\line(3,2){24.2}}
\put(215.3,295){\circle*{5}} \put(230,250){\line(1,3){9}}
\put(239,278.5){\circle*{5}} \put(239,278.5){\line(-3,2){24.2}}
\put(240,250){,} \put(250,250){$\ldots$}
\end{picture}
\caption{Extended Dynkin diagrams $\widetilde{A}_n$,
$\widetilde{D}_n$, $\widetilde{E}_6$, $\widetilde{E}_7$,
$\widetilde{E}_8$.}\label{Ext_Dynk_diagr}
\end{figure}

We also need the following statement.
\begin{theorem}[J.H.\ Smith]\label{theorem:Smith}
Let $\Gamma$ be a graph with index $r$. Then $r\leq 2$ $(r<2)$ if
and only if each component of~$\Gamma$ is a subgraph (proper
subgraph) of one of the graphs depicted in Fig.~{\rm
\ref{Ext_Dynk_diagr}} which all have an index equal to~$2$.
\end{theorem}

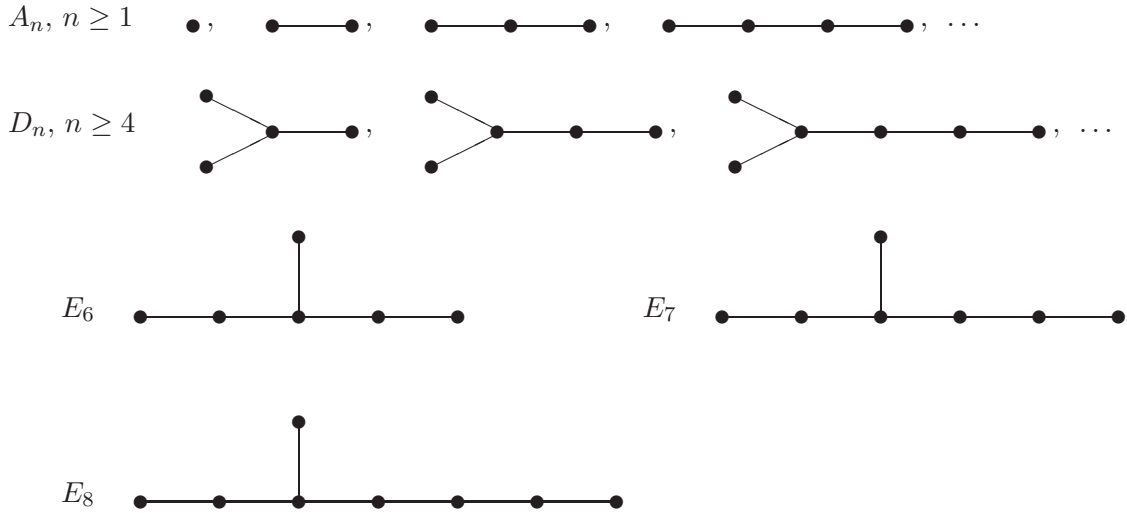
\begin{figure}[t]
\hspace*{5mm}\begin{picture}(440,200) \put(0,180){$A_n,\,n\geq 1$}
\put(70,180){\circle*{5}} \put(75,180){,}
\put(100,180){\circle*{5}} \put(100,180){\line(1,0){30}}
\put(130,180){\circle*{5}} \put(135,180){,}
\put(160,180){\circle*{5}} \put(160,180){\line(1,0){30}}
\put(190,180){\circle*{5}} \put(190,180){\line(1,0){30}}
\put(220,180){\circle*{5}} \put(225,180){,}
\put(250,180){\circle*{5}} \put(250,180){\line(1,0){30}}
\put(280,180){\circle*{5}} \put(280,180){\line(1,0){30}}
\put(310,180){\circle*{5}} \put(310,180){\line(1,0){30}}
\put(340,180){\circle*{5}} \put(345,180){,}
\put(355,180){$\ldots$}
\put(0,140){$D_n,\, n\geq 4$} \put(75,153.41){\circle*{5}}
\put(75,153.41){\line(2,-1){26.8}} \put(75,126.59){\circle*{5}}
\put(75,126.59){\line(2,1){26.8}} \put(100,140){\circle*{5}}
\put(100,140){\line(1,0){30}} \put(130,140){\circle*{5}}
\put(135,140){,}  
\put(160,153.41){\circle*{5}} \put(160,153.41){\line(2,-1){26.8}}
\put(160,126.59){\circle*{5}} \put(160,126.59){\line(2,1){26.8}}
\put(185,140){\circle*{5}} \put(185,140){\line(1,0){30}}
\put(215,140){\circle*{5}} \put(215,140){\line(1,0){30}}
\put(245,140){\circle*{5}} \put(250,140){,} 
\put(275,153.41){\circle*{5}} \put(275,153.41){\line(2,-1){26.8}}
\put(275,126.59){\circle*{5}} \put(275,126.59){\line(2,1){26.8}}
\put(300,140){\circle*{5}} \put(300,140){\line(1,0){30}}
\put(330,140){\circle*{5}} \put(330,140){\line(1,0){30}}
\put(360,140){\circle*{5}} \put(360,140){\line(1,0){30}}
\put(390,140){\circle*{5}} \put(395,140){,}
\put(405,140){$\ldots$}
\put(20,70){$E_6$} \put(50,70){\circle*{5}}
\put(50,70){\line(1,0){30}} \put(80,70){\circle*{5}}
\put(80,70){\line(1,0){30}} \put(110,70){\circle*{5}}
\put(110,70){\line(0,1){30}} \put(110,100){\circle*{5}}
\put(110,70){\line(1,0){30}} \put(140,70){\circle*{5}}
\put(140,70){\line(1,0){30}} \put(170,70){\circle*{5}}
\put(240,70){$E_7$} \put(270,70){\circle*{5}}
\put(270,70){\line(1,0){30}} \put(300,70){\circle*{5}}
\put(300,70){\line(1,0){30}} \put(330,70){\circle*{5}}
\put(330,70){\line(0,1){30}} \put(330,100){\circle*{5}}
\put(330,70){\line(1,0){30}} \put(360,70){\circle*{5}}
\put(360,70){\line(1,0){30}} \put(390,70){\circle*{5}}
\put(390,70){\line(1,0){30}} \put(420,70){\circle*{5}}
\put(20,0){$E_8$} \put(50,0){\circle*{5}}
\put(50,0){\line(1,0){30}} \put(80,0){\circle*{5}}
\put(80,0){\line(1,0){30}} \put(110,0){\circle*{5}}
\put(110,0){\line(0,1){30}} \put(110,30){\circle*{5}}
\put(110,0){\line(1,0){30}} \put(140,0){\circle*{5}}
\put(140,0){\line(1,0){30}} \put(170,0){\circle*{5}}
\put(170,0){\line(1,0){30}} \put(200,0){\circle*{5}}
\put(200,0){\line(1,0){30}} \put(230,0){\circle*{5}}
\end{picture}
\caption{Dynkin diagrams $A_n$, $D_n$, $E_6$, $E_7$,
$E_8$.}\label{Dynk_diagr}
\end{figure}

\begin{corollary}[For trees]\label{trees}
Let $\Gamma$ be a tree with index $r.$ Then

$1.$  $r<2$ if and only if $\Gamma$ is one of the following
graphs: $A_n$, $D_n$, $E_6$, $E_7$, $E_8$ (see Fig.~{\rm
\ref{Dynk_diagr}}).

$2.$ $r=2$ if and only if $\Gamma$ is one of the following graphs:
$\widetilde{D}_n$, $\widetilde{E}_6$, $\widetilde{E}_7$,
$\widetilde{E}_8 $ (see Fig.~{\rm \ref{Ext_Dynk_diagr}}).
\end{corollary}

\section[Configurations of subspaces in a Hilbert space with fixed angles between
them]{Conf\/igurations of subspaces in a Hilbert space\\ with
f\/ixed angles between them}\label{Subsp_in_Hilbert_spase}

Let $H$ be a complex Hilbert space and let $H_i,H_j\subset H$ be
its closed subspaces. We say that an angle between $H_i$ and $H_j$
is f\/ixed and equals to $\varphi_{i,j}\in [0;\frac{\pi}{2}]$ if
for the orthogonal projections $P_{H_i}$, $P_{H_j}$ on these
subspaces  we have
\[
P_{H_i}P_{H_j}P_{H_i}=\cos^2(\varphi_{i,j}) P_{H_i} \qquad
\text{and}\qquad P_{H_j}P_{H_i}P_{H_j}=\cos^2(\varphi_{i,j})
P_{H_j}.
\]
Having a f\/inite undirected graph $\Gamma$ without multiple edges
and loops with the numbers on its edges, we def\/ine the
conditions on the conf\/iguration of subspaces in a Hilbert space
as follows. The subspaces correspond to the vertices of a graph
and an angle between any two of them is given by the number
$\tau_{i,j}$ standing on the respective edge. If vertices are not
adjacent we assume that the corresponding subspaces are
orthogonal.

We consider the following questions:
\begin{enumerate}\itemsep=0pt\vspace{-2mm}
\item For which values of the parameters $\tau_{i,j}$ the
conf\/iguration associated with a graph $\Gamma$ exists. \item
Give the description of all irreducible conf\/igurations
(associated with a f\/ixed graph $\Gamma$ and an arrangement of
numbers on its edges) up to a unitary transformation.
\end{enumerate}

It should be noted that the subspaces corresponding to vertices
from dif\/ferent components of~$\Gamma$ are orthogonal, so we will
consider only connected graphs.

These problems can be reformulated in terms of f\/inding
$*$-representations of $*$-algebras associated with $\Gamma$ with
the numbers $\tau_{i,j}$ on the edges. Let $\Gamma$ be a f\/inite,
undirected, connected graph without multiple edges and loops, with
$\Gamma_0$ $(|\Gamma_0|=n)$ and $\Gamma_1$ the sets of the
vertices and the edges respectively. Let $\tau :\Gamma_1
\rightarrow (0,1)$ be the arrangement of numbers on its edges. We
enumerate the vertices of $\Gamma$ by numbers $1,\ldots, n$ in any
way and denote $\tau (i,j)=:\tau_{i,j}=\tau_{j,i}.$
\begin{definition}
$A_{\Gamma,\tau}$ is an $*$-algebra with \textbf{1} over
$\mathbb{C}$ generated by projections $p_1,\ldots , p_n$
($p_i^2=p_i^*=p_i$, $\forall\, i$) with relations
\begin{gather*}
p_i p_jp_i=\tau_{i,j}p_i \quad \text{and}\quad
p_jp_ip_j=\tau_{i,j}p_j \ \ \text{if} \ \  (i,j)\in\Gamma_1, \quad
\text{and}\quad p_i p_j=p_jp_i=0 \ \ \text{otherwise}.
\end{gather*}
\end{definition}

Results on dimension of the algebra $A_{\Gamma,\tau}$ (in
dependence on a graph $\Gamma$) can be found in
\cite{MashaVlasenko,Vlasenko&Popova}.

Below we suppose that $\Gamma$ is a \textit{tree}. Then the
$*$-algebra $A_{\Gamma,\tau}$ is f\/inite dimensional and,
therefore, does not have inf\/inite dimensional irreducible
$*$-representations.

\section[On the set of values of the parameters where $A_{\Gamma,\tau}$
has  $*$-representations]{On the set of values of the parameters\\
where $\boldsymbol{A_{\Gamma,\tau}}$ has
$\boldsymbol{*}$-representations}\label{OnSetSigma}

Let $\Gamma$ be a tree and $\mathcal{A}(\Gamma,\tau)=\|
\mathcal{A}_{i,j}\|_{i,j=1}^n$ be the symmetric matrix with
$\mathcal{A}_{i,i}=1$ $ \forall\, i$;
$\mathcal{A}_{i,j}=\sqrt{\tau_{i,j}}$ if $(i,j)\in\Gamma_1$, and
$\mathcal{A}_{i,j}=0$ otherwise.

\begin{proposition}
Let $\Gamma$ be a tree. Nontrivial $*$-representations of an
algebra $A_{\Gamma,\tau}$ exist if and only if the matrix
$\mathcal{A}(\Gamma,\tau)$ is positive semidefinite. Irreducible
nontrivial $*$-representation is unique up to the unitary
equivalence and its dimension is equal to the
$\mathrm{rank}\,(\mathcal{A}(\Gamma,\tau)).$
\end{proposition}
The proof one can f\/ind in \cite{Vlasenko&Popova}.

In the following we suppose that $\tau$ is \textit{constant}
($\tau_{i,j}=\tau$ $\forall \, (i,j)\in\Gamma_1$). The set of
values of the parameter $\tau$ where $A_{\Gamma,\tau}$ has
nontrivial $*$-representations we will denote by
$\Sigma_{\Gamma}.$

\begin{theorem}\label{Theorem:SigmaG}
Let $\Gamma$ be a tree with index $r.$ Then
\begin{gather}\label{eq3}
\Sigma_{\Gamma}=\left(0;\frac{1}{r^2}\right].
\end{gather}
\end{theorem}
\begin{proof}
Indeed, $\mathcal{A}(\Gamma,\tau)=I+\sqrt{\tau}B_{\Gamma}$ where
$I$ is $n\times n$ unit matrix and $B_{\Gamma}$ is the adjacency
matrix of the tree $\Gamma.$ The matrix $\mathcal{A}(\Gamma,\tau)$
is positive semidef\/inite if and only if its minimal eigenvalue
is nonnegative, i.e., $1+\sqrt{\tau}q\geq 0$ (in the notations of
Section \ref{about_graphs}) which is equivalent to
$\tau\leq\frac{1}{q^2}.$
 For trees we know that $q=-r$ (see Remark \ref{q=-r}), so
the theorem is proved.
\end{proof}

\begin{example}
Let us f\/ind the sets $\Sigma_\Gamma$ where graphs $\Gamma$ are
Dynkin diagrams.
\begin{gather*}
\Sigma_{A_n}=\left(0;\frac{1}{4\cos^2
\frac{\pi}{n+1}}\right],\qquad
\Sigma_{D_n}=\left(0;\frac{1}{4\cos^2 \frac{\pi}{2(n-1)}}\right],\\
\Sigma_{E_6}=\left(0;\frac{1}{4\cos^2
\frac{\pi}{12}}\right],\qquad
\Sigma_{E_7}=\left(0;\frac{1}{4\cos^2
\frac{\pi}{18}}\right],\qquad
\Sigma_{E_8}=\left(0;\frac{1}{4\cos^2 \frac{\pi}{30}}\right].
\end{gather*}
Values of index of Dynkin diagrams can be found
in~\cite{cvetkovic}.
\end{example}

Now we give some properties of $\Sigma_{\Gamma}$ ($\Gamma$ is a
tree) that immediately  follow from Proposition~\ref{inequalities}
and Theorem~\ref{Theorem:SigmaG}.

\begin{proposition}
Let $\Gamma$ be a tree with $n$ vertices. Then
\[
1) \ \ \left(0;\frac{1}{(n-1)^2}\right]\subseteq \Sigma_{\Gamma},
\qquad 2) \ \ \Sigma_{\Gamma}\subseteq\left(0;\frac{1}{4\cos^2
\frac{\pi}{n+1}}\right].
\]
\end{proposition}
\begin{remark}\label{rem2}
For graph $\Gamma$ that is not a tree situation is more
complicated. For example, if the graph is a cycle with $n$
vertices, i.e.\ $\Gamma=\widetilde{A}_{n-1},$ we know that
$\Sigma_{\Gamma}=\Sigma_{A_{n-1}}=\left(0;\frac{1}{4\cos^2
\frac{\pi}{n}}\right]$ (see~\cite{Popova}). But the index of
$\widetilde{A}_{n-1}$ is $r=2$ (see Theorem~\ref{theorem:Smith})
and formula \eqref{eq3} does not hold. Moreover, it is known that
all eigenvalues of $\widetilde{A}_{n-1}$ are of the form:
$\lambda_j=2\cos \frac{2\pi}{n}j$, $j=1,\ldots ,n$
(see~\cite{cvetkovic}). Therefore, if n is even then no one
eigenvalue $\lambda_j$ of the graph $\widetilde{A}_{n-1}$ makes
the formula $\Sigma_{\Gamma}=\Big(0;\frac{1}{\lambda_j^2}\Big]$
true.
\end{remark}

Next proposition follows directly from Corollary~\ref{trees} and
Theorem~\ref{Theorem:SigmaG}.

\begin{proposition}
Let $\Gamma$ be a tree. Then
\begin{enumerate}\itemsep=0pt
\item[{\rm 1.}] $\max\Sigma_{\Gamma}> 1/4$ if and only if $\Gamma$
is one of the following graphs: $ A_n$, $D_n$, $E_6$, $E_7$,
$E_8.$ \item[{\rm 2.}] $\max\Sigma_{\Gamma}=1/4$ if and only if
$\Gamma$ is one of the following graphs: $\widetilde{D}_n$,
$\widetilde{E}_6$, $\widetilde{E}_7$, $\widetilde{E}_8.$
\item[{\rm 3.}] For all other trees which are neither Dynkin
diagrams nor extended Dynkin diagrams we have
$\max\Sigma_{\Gamma}< 1/4.$
\end{enumerate}
\end{proposition}

\LastPageEnding

\end{document}